\documentclass[11pt]{article}
\usepackage{}

\usepackage[paperwidth=10.5in,top=1.3in, bottom=1.2in]{geometry}
\usepackage{amsfonts}
\usepackage{amsthm}
\usepackage{amssymb}
\usepackage{amsmath,amsfonts}
\usepackage{mathrsfs}
\usepackage{cases}
\usepackage{latexsym,bm}
\usepackage{indentfirst}
\usepackage{color}
\usepackage{ifpdf}
\usepackage{graphicx}
\usepackage{psfrag}
\usepackage{dsfont}
\usepackage[
pdfauthor={CSY},
pdftitle={Strong generalized Halin graphs},
pdfstartview=XYZ,
bookmarks=true,
colorlinks=true,
linkcolor=blue,
urlcolor=blue,
citecolor=blue,
bookmarks=true,
linktocpage=true,
hyperindex=true
]{hyperref}

\newcommand{\xiaowuhao}{\fontsize{9pt}{\baselineskip}\selectfont}

\newtheorem{THM}{\textbf{Theorem}}[section]

\newtheorem{LEM}{\textbf{Lemma}}[section]
\newtheorem{CON}{\textbf{Conjecture}}

\newcommand{\pf}{\textbf{Proof}.\quad}
\newcommand{\spf}{ \emph{Proof}.\quad}

\newtheorem{CLA}{\textbf{Claim}}[section]
\newcommand{\qqed}{\hfill $\blacksquare$\vspace{1mm}}

\newcommand{\CC}{\mathcal{C}}
\newcommand{\HH}{\mathcal{H}}

\begin{document}

\title{Vizing's 2-factor Conjecture Involving Large Maximum Degree }
\author{ Guantao Chen and Songling Shan\\
{\xiaowuhao  Georgia State University, Atlanta, GA\,30303, USA}}
\date{}
\maketitle
\emph{\textbf{Abstract}.}
Let $G$ be a connected simple graph of order $n$ and let $\Delta(G)$
and $\chi'(G)$ denote the maximum degree and chromatic index of
$G$, respectively. Vizing proved that $\chi'(G)=\Delta(G)$ or $\Delta(G)+1$.
Following this result, $G$ is called $\Delta$-critical if $\chi'(G)=\Delta(G)+1$ and $\chi'(G-e)=\Delta(G)$ for every $e\in E(G)$. In 1968, Vizing conjectured that if $G$ is an $n$-vertex
$\Delta$-critical graph, then the
independence number $\alpha(G)\le n/2$. Furthermore, he conjectured that, in fact, $G$ has a 2-factor.
Luo and Zhao showed that if $G$ is an $n$-vertex $\Delta$-critical
 graph with  $\Delta(G)\ge n/2$, then $\alpha(G)\le n/2$. More recently,
they showed that if $G$ is an $n$-vertex $\Delta$-critical graph with
$\Delta(G)\ge
6n/7$,  then $G$ has  a hamiltonian cycle, and so $G$ has a 2-factor.
In this paper, we show that if $G$ is an $n$-vertex $\Delta$-critical graph with $\Delta(G)\ge n/2$,
then $G$ has a 2-factor.


\vspace{2mm}

\section{Introduction}

In this paper, we only consider simple and finite graphs. Let
$G$ be a graph. We fix the notation $\Delta$ for the maximum degree of $G$ throughout the paper.
A $k$-vertex of $G$ is a vertex of degree
 $k$ in $G$. Denote by  $V_{\Delta}$ the set of $\Delta$-vertices in $G$ and by $\chi'(G)$ the edge-chromatic index of   $G$.
The graph $G$ is called {\it critical}\,({\it edge-chromatic critical}) if it has no isolated vertices and
$\chi'(G-e)<\chi'(G)$ for every  $e\in E(G)$. From the definition, it is clear that if $G$ is critical,
then $G$ is connected.
In 1965, Vizing~\cite{Vizing-2-classes} showed that a graph of maximum degree
$\Delta$ has edge chromatic index either $\Delta$ or $\Delta+1$.
If $\chi'(G)=\Delta$, then $G$ is said to be of class 1; otherwise, it is said to be
of class 2.  Appearing  easily, however,
Holyer~\cite{Holyer} showed that it is NP-complete to determine whether an arbitrary graph is of class 1.
A critical graph $G$ is called {\it $\Delta$-critical} if
$\chi'(G)=\Delta+1$.
So $\Delta$-critical graphs are
class 2 graphs. On the other hand,  every critical  class 2 graph of maximum degree $\Delta$ is  a
$\Delta$-critical graph.
Motivated by the classification problem, Vizing studied critical class 2 graphs and made the following two well-known conjectures.

The first one, appeared in~\cite{vizing-ind},
is on  the independence number $\alpha(G)$ of $G$, that is, the size of a maximum independent set
of $G$.

\begin{CON}[Vizing's Independence Number Conjecture]\label{ind}
Let $G$ be a $\Delta$-critical graph of order $n$.
Then $\alpha(G)\le n/2$.
\end{CON}

The second one, appeared in~\cite{vizing-2factor}, is on 2-factor, a 2-regular spanning
subgraph.

\begin{CON}[Vizing's 2-factor Conjecture]\label{2-factor}
Let $G$ be a $\Delta$-critical graph. Then $G$ contains a 2-factor.
\end{CON}

%

%
%
As each cycle $C$ satisfying $\alpha(C)\le |V(C)|/2$,  Conjecture~\ref{2-factor} implies Conjecture~\ref{ind}.
For the Independence Number Conjecture,
Brinkmann et al.~\cite{ind1} in 2000
 proved that if G   is an $n$-vertex  $\Delta$-critical graph, then $\alpha(G)<2n/3$; and
 the upper bound is further improved when
 the maximum degree is between 3 and 10. In 2006,
 Luo and Zhao~\cite{vizing-independence-large-Delta} confirmed the conjecture for graphs
 with large maximum degree.

\begin{THM}\label{independence}
Let $G$ be an $n$-vertex $\Delta$-critical graph.
 Then $\alpha(G)\le n/2$ if $\Delta\ge n/2$.
\end{THM}

Additionally, Luo and Zhao~\cite{ind2} in 2008
showed that if $G$ is an $n$-vertex
$\Delta$-critical graph, then
$\alpha(G)<(5\Delta-6)n/(8\Delta-6)<5n/8$ when $\Delta\ge6$.
In 2009,
Woodall~\cite{ind3} further improved the upper bound to
$3n/5$.  Compared to the progress on  Vizing's Independence Number Conjecture,
the progress on the 2-factor Conjecture is slower.
In 2004,
Gr{\"u}newald and Steffen~\cite{vizing-2-factor-overful} established Vizing's
2-factor conjecture for graphs with the deficiency $\sum_{v\in V(G)}(\Delta(G)-d_G(v))$
small;
in particular, for overfull graphs, i.e., graphs
of odd order and with the deficiency $\sum_{v\in V(G)}(\Delta(G)-d_G(v))<\Delta(G)$.
In 2012, Luo and Zhao~\cite{Vizing-2-factor-hamiltonian}
proved that if $G$ is an $n$-vertex $\Delta$-critical graph with $\Delta\ge 6n/7$,
then $G$ contains a Hamiltonian cycle, and thus a 2-factor with exactly one component.
Still considering $\Delta$-critical graphs with large maximum degree,
in line with Luo and Zhao's result on the Independence Number Conjecture\,(Theorem~\ref{independence}),
in this paper, we reduce the lower bound from $6n/7$ to $n/2$ as follows.

%
%
%

\begin{THM}\label{main}
Let $G$ be an $n$-vertex  $\Delta$-critical graph. Then $G$
has a 2-factor if $\Delta\ge n/2$.
\end{THM}

\section{Notations and Lemmas}

For a vertex $x$ of a graph $G$,
we denote by $N_G(x)$ the set of neighbors of $x$ in $G$ and by $d_G(x)$ the degree of $x$ in $G$.
For a set of vertices $S$ in $G$,
we define $N_G(S)$ by
$N_G(S)=\bigcup_{x\in S} N_G(x)$.
For disjoint sets of vertices $S$ and $T$ in $G$,
we denote by $e_G(S, T)=|E_G(S,T)|$, the number of edges that has one
end vertex in $S$ and the other in $T$.
If $S$ is a singleton set $S=\{s\}$,
we write $e_G(s, T)$ instead of $e_G(\{s\}, T)$.
If $G$ is a bipartite graph with partite sets $A$ and $B$, we
denote $G$ by $G[A,B]$ to emphasize the two partite sets.
To prove Theorem~\ref{main}, we present a few lemmas.

\begin{LEM}[Vizing's Adjacency Lemma]\label{vizing adjacency lemma}
Let $G$ be a $\Delta$-critical graph. Then for any edge
$xy\in E(G)$, $x$ is adjacent to at least $\Delta-d_G(y)+1$  $\Delta$-vertices
$z$ with $z\ne y$.
\end{LEM}

%

As there are two specified  bipartite graphs $H^*[X',T]$ and
$H[X,T]$ defined in the sequel, for consistency,  we use notation $H^*[X',T]$
in  lemmas only regarding to the bipartite graph $H^*[X',T]$.
Denote by $R[A,B]$ for a general bipartite graph in distinguishing
with the  bipartite graphs $H^*$ and $H$.
A matching of a graph $G$ is a set of independent
edges in $G$.
If $M$ is a matching of $G$,
then let $V(M)$ denote the set of end vertices
of the edges in $M$.
For $X\subseteq V(G)$,
$M$ is said to saturate $X$ if $X\subseteq V(M)$.
The following result,  which guarantees a matching
saturating  one partite set in a bipartite graph,
can also be found in~\cite{local-CE}.



\begin{LEM}\label{matching}
Let $H^*$ be a bipartite graph with partite sets
$X'$ and $T$.
If there is no isolated vertex in $T$
and $d_{H^*}(y)\ge d_{H^*}(x)$
for every edge $xy$ with $x\in X'$ and $y\in T$,
then $H^*$ has a matching which saturates $T$.
\end{LEM}

\pf Suppose not. Then by Hall's Theorem, there
is a nonempty set $A\subseteq T$ such that
$|N_{H^*}(A)|<|A|$. We choose  $A$ such that
it has the minimum cardinality under the constraint that
$|N_{H^*}(A)|<|A|$. Let $B:=N_{H^*}(A)$ and $H':=H^*[A\cup B]$ be
the subgraph induced by $A\cup B$.
We claim that, in $H'$,  there is a matching saturating
$B$. Suppose not. Then
by Hall's Theorem again, there is a nonempty subset $B'\subseteq B$
such that $|N_{H'}(B')|<|B'|$.
Since $B'\subseteq B=N_{H^*}(A)\ne \emptyset$\,($T$ has no isolated vertices), $N_{H'}(B')\ne \emptyset$.  Let $A'=A-N_{H'}(B')$.
As $|A|>|B|\ge |B'|>|N_{H'}(B')|>0$, we have $0<|A'|<|A|$.
On the other hand, we have $N_{H'}(A')=N_{H^*}(A')=B-B'$.
So, the sequence of inequalities  $|A'|=|A|-|N_{H'}(B')|>|B|-|N_{H'}(B')|> |B|-|B'|=|B-B'|=|N_{H^*}(A')|$ holds,
showing a contradiction to the minimality of $A$ under the condition $|N_{H^*}(A)|<|A|$.
Let $M$ be a matching which saturates $B$ in $H'$.
Since $|A|>|B|$, $A-V(M)\ne \emptyset$.
Let $y^*\in A-V(M)$. Then as $T$ has no isolated vertices,
$d_{H^*}(y^*)=d_{H'}(y^*)\ge 1$. Thus,
\begin{eqnarray*}
 e_{H^*}(A,B) &=& \sum_{xy\in M, x\in B, y\in A}d_{H'}(x)\quad\mbox{($M$ saturates $B$ in $H'$)}    \\
 &\le & \sum_{xy\in M, x\in B, y\in A}d_{H^*}(x)\\
   &\le &\sum_{xy\in M, x\in B, y\in A}d_{H^*}(y)  \\
   &<& \sum_{xy\in M, x\in B, y\in A}d_{H^*}(y)+d_{H^*}(y^*)\\
   &\le&  e_{H^*}(A,B),
\end{eqnarray*}
showing a contradiction.
\qqed

%

\begin{LEM}\label{islatedmatching}
Let $R$ be a bipartite graph with
partite sets $A$ and $B$,
$A_1=\{x\in A\,|\, d_R(x)=1\}$ and $B_1=N_R(A_1)$.
Then $R$ has a matching saturating $B$ if the bipartite graph
$R'[A-A_1,B-B_1]:=R[(A-A_1)\cup (B-B_1)]$ has a matching  saturating $B-B_1$.
\end{LEM}

\pf Suppose that $R'$ has a matching $M'$ which saturates $B-B_1$.
Since each vertex in $A_1$ has a unique neighbor in $B_1$,
there is a matching $M_0$ saturating $B_1$ in the subgraph of
$R$ induced on $A_1\cup B_1$.
Then $M'\cup M_0$ is a matching which saturates $B$ in $R$.
\qqed

The following lemma
is a generalization of a
result in~\cite{vizing-independence-large-Delta}.

\begin{LEM}\label{TnoDelta}
Let $G$ be a $\Delta$-critical graph and
$T$ be an independent set of $G$.
Let $X'=V(G)-T$ and $H^*:=G-E(G[X'])$ be the bipartite
graph with partite sets $X'$ and $T$.
Then for each edge
$xy\in E(H^*)$ with $x\in X'$ and $y\in T$,
$d_{H^*}(y)\ge d_{H^*}(x)+1-\delta_0+\sigma_x$,
where $\delta_0=|T\cap V_{\Delta}|$
is the number of $\Delta$-degree vertices in $T$ and
$\sigma_x$ is the number of non $\Delta$-degree neighbors of $x$ in $X'$.
Moreover, if $\delta_0\le 1$, then
there is a matching which saturates $T$ in $H^*$.
\end{LEM}
\pf
Let $xy\in E(H^*)$ with $x\in X'$ and $y\in T$.
By Vizing's Adjacency Lemma\,(Lemma~\ref{vizing adjacency lemma}),
$x$ is adjacent to at least $\Delta-d_G(y)+1$  $\Delta$-vertices
in $G$.
As $T$ has  $\delta_0$ $\Delta$-vertices,
we know $x$ is adjacent to
at least $\Delta-d_G(y)+1-\delta_0$ $\Delta$-vertices in
$X'$\,(notice that the quantity is meaningful only if
$\Delta-d_G(y)+1-\delta_0>0$).
Then, $\Delta\ge d_G(x)=d_{H^*}(x)+e_G(x,X')\ge d_{H^*}(x)+\Delta-d_G(y)+1-\delta_0+\sigma_x $.
Thus $d_{H^*}(y)=d_G(y)\ge d_{H^*}(x)+1-\delta_0+\sigma_x$.

When $\delta_0\le 1$,  for every edge
$xy\in E(H^*)$ with $x\in X'$ and $y\in T$,
the inequalities $d_{H^*}(y)\ge d_{H^*}(x)+\sigma_x\ge d_{H^*}(x)$
hold. As $G$ is $\Delta$-critical, it is connected.
Consequently, $d_{H^*}(y)=d_G(y)\ge 1$ for each $y\in T$.
By applying Lemma~\ref{matching},
we see that there is a matching which saturates $T$ in
$H^*$.
\qqed

\section{A Detour to Tutte's 2-Factor Theorem}

Tutte in~\cite{tutte-2-factor} obtained a necessary and
sufficient condition for a graph to contain an $f$-factor;
the characterization involves pairs of
two disjoint vertex sets. Let $G$ be a graph and
$(S,T)$ be an ordered pair of disjoint vertex sets of $G$.
A  component $C$ of $G-(S\cup T)$
is said to be an {\it odd component} w.r.t. $(S,T)$
(resp.~{\it even component} w.r.t. $(S,T)$)
if $e_G(C, T)\equiv 1\pmod{2}$
(resp.~$e_G(C, T)\equiv 0\pmod{2}$).
Let $\HH_G(S, T)$ be the set of odd components of $G-(S\cup T)$,
$h_G(S, T)=|\HH_G(S, T)|$,
and let
$\delta_G(S, T)=2|S|+\sum_{v\in T} d_{G-S}(v)-2|T|-h_G(S, T)$.
It is easy to see $\delta_G(S, T)\equiv 0\pmod{2}$
for every $S$,~$T\subseteq V(G)$
with $S\cap T=\emptyset$.
We use the following criterion for the existence of a $2$-factor,
which is a restricted form of Tutte's $f$-Factor Theorem.

\begin{THM}\label{tutte's theorem}
A graph $G$ has a $2$-factor if and only if
$\delta_G(S, T)\ge 0$
for every $S$,~$T\subseteq V(G)$
with $S\cap T=\emptyset$.
\end{THM}

%

An ordered pair $(S,T)$ consists of
 disjoint sets  of vertices $S$ and $T$ in a graph $G$
is called a barrier if $\delta_G(S, T)\le -2$.
By Theorem~\ref{tutte's theorem},
every graph $G$ without a 2-factor has a barrier.
A barrier $(S, T)$ is called a minimum barrier
if $|S\cup T|$ is  smallest among all the barriers of $G$.
A minimum barrier
of a graph without a $2$-factor has some nice properties,
see~\cite{AEFOS,local-CE} for examples.
We will use the  properties listed  in the
following lemma in our proof.

\begin{LEM}\label{minimal_barrier}
Let $G$ be a graph without  a $2$-factor
and $(S, T)$ be a minimum barrier of $G$.
Then the following statements hold.
\begin{enumerate}
\item
$T$ is independent,
\item for every even component $C$ w.r.t. $(S,T)$, $e_G(T, V(C))=0$, and
\item
for every odd component $C$ w.r.t. $(S,T)$
and every $v\in T$, $e_G(v, V(C)) \le 1$, i.e.,
either $e_G(v, V(C))=0$ or $e_G(v, V(C)) = 1$.
\end{enumerate}
\end{LEM}

Let $(S,T)$ be a minimum barrier of $G$. We introduce some necessary notations
w.r.t. $(S,T)$ for this paper. Denote
$$
\CC_k=\{C\in \HH_G(S, T)\,|\, e_G(T,V(C))=k\}.
$$
Then $\HH_G(S, T)=\bigcup_{k\ge 0}\CC_{2k+1}$ and $h_G(S,T)=|\cup_{k\ge 0}\CC_{2k+1}|$.
For any  $v\in T$, let
$$
\CC_v  = \{C\in \HH_G(S, T)\,|\, e_G(v,V(C))=1\} \quad \text{and}\quad  \CC_{1v}  =  \{C\in \CC_1 \,|\, e_G(v,V(C))=1\}.
$$
It is clear that $\CC_{1v}\subseteq \CC_{v}$. We distinguish $\CC_{1v}$ because in
the proof of Theorem~\ref{main}, we need pay special attention on
vertices $v\in T$ with $\CC_{1v}\ne \emptyset$.


\begin{LEM}\label{size2_0}
Let $G$ be a graph without a 2-factor and $(S,T)$ be a minimum barrier
such that $h_G(S,T)$ is smallest,
and let $v\in T$ with $|\CC_v|\ge 2$ and $|\CC_{1v}|\ge 1$.
Then for any vertex $w$ in a component $D\in \CC_{1v}$,
$e_G(w,V(D)\cup \{v\})\ge 2$.
\end{LEM}

\spf
Suppose on the contrary that 
there exists  a vertex $w\in V(D)$
 such that $e_G(w, V(D)\cup \{v\})\le 1$.  Since
 $G^*:=G[V(D)\cup \{v\}]$ is connected, $e_G(w, V(D)\cup \{v\})=1$,
which in turn gives that $G^*-w$ is connected.
 Let $T^*=(T-\{v\})\cup \{w\} $.
We claim that $(S,T^*)$ is a minimum barrier with
$h_G(S,T^*)<h_G(S,T)$. This will give a contradiction to
the choice of $(S,T)$.

To see that  $(S,T^*)$ is a minimum barrier we calculate
$\delta_G(S,T^*)$.
Let
$D_{v}$ be the component of $G-(S\cup T^*)$
containing $v$. Notice that besides $v$, the component $D_v$
contains also vertices in $D-w$ and
all the $d_{G-S}(v)-1$ odd components $C\,(\ne D)$ of $G-(S\cup T)$
 such  that $e_G(v, V(C))=1$.

We first show that $D_v$ is an odd component of $G-(S\cup T^*)$.
Let $v^*$ be the neighbor of $v$ in $D$.
For each $C\in\HH_G(S,T)\cap \CC_v$, denote
the odd number $e_G(T,V(C))$ by $2k_c+1$
for some nonnegative integer $k_c$.
If $w=v^*$, then as $e_G(w, V(D)\cup \{v\})=1$,
we see that $D$ is a single vertex component and  $V(G^*)=\{v,w\}$. Then
\begin{eqnarray*}
 e_G(T^*,V(D_{v})) &=& e_G(T,V(D_v))-(d_{G-S}(v)-1)+e_G(v^*,v) \\
 &=&\sum\limits_{C\in \CC_v-\{D\} }e_G(T,V(C))-(d_{G-S}(v)-1)+e_G(v^*,v)\\
   &=& \sum\limits_{C\in \CC_v-\{D\}}(2k_c+1-1)+e_G(v^*,v)\quad \mbox{(noticing that $|\CC_v-\{D\}|=d_{G-S}(v)-1$)},
\end{eqnarray*}
which is odd by $e(v^*,v)=d_{G-S}(v^*)=1$. So
$D_{v}$ is an odd component of
$G-(S\cup T^*)$.
If $w\ne v^*$, then $D$ has at least two vertices, and
\begin{eqnarray*}
 e_G(T^*,V(D_{v})) &=& e_G(T,V(D_v))-d_{G-S}(v)+e_G(w,V(G^*)) \\
 &=&\sum\limits_{C\in \CC_v }e_G(T,V(C))-d_{G-S}(v)+e_G(w,V(G^*))\\
   &=& \sum\limits_{C\in \CC_v}(2k_c+1-1)+1\quad \mbox{(noticing that $|\CC_v|=d_{G-S}(v)$)},
\end{eqnarray*}
which is again odd.
Hence,
$h_G(S,T^*)=h_G(S,T)-d_{G-S}(v)+1$.
So,
\[
\begin{split}
\delta_G(S, T^*)&=2|S|-2|T^*|+\sum_{y\in T^*}d_{G-S}(y)-h_G(S, T^*)\\
&=2|S|-2|T|+\sum_{y\in T}d_{G-S}(y)-d_{G-S}(v)+d_{G-S}(w)-
(h_G(S,T)-d_{G-S}(v)+1)\\
&=2|S|-2|T|+\sum_{y\in T}d_{G-S}(y)-h_G(S,T)\\
&\le -2;
\end{split}
\]
by noticing that $d_{G-S}(w)=d_{G^*}(w)=1$.

As $|S\cup T^*|=|S\cup T|$, $(S,T^*)$ is a minimum barrier.
However, as $d_{G-S}(v)= |\CC_{v}|\ge 2$,
 $h_G(S,T^*)=h_G(S,T)-d_{G-S}(v)+1<h_G(S,T)$.
 This gives a contradiction to the choice of $(S,T)$.
\qqed

The following result is a consequence  of Lemma~\ref{size2_0}.

\begin{LEM}\label{size2}
Let $G$ be a graph without a 2-factor and $(S,T)$ be a minimum barrier
such that $h_G(S,T)$ is smallest.
Then for any  $v\in T$ with $|\CC_{v}|\ge 2$ and
$D\in \CC_{1v}$\,(if exists),
$|V(D)|\ge 2$.
\end{LEM}



\section{Proof of the Main Result}

Assume,
to the contrary,
that the $n$-vertex $\Delta$-critical graph $G$ with $\Delta \ge n/2$ does not have a $2$-factor.
Then $\Delta\ge 3$ since a 2-critical graph is an odd cycle, which is a 2-factor of $G$.
Since $G$ is $\Delta$-critical, by Vizing's Adjacency Lemma,
each vertex of $G$ is adjacent to at least two $\Delta$-vertices and thus $\delta(G)\ge 2$.

By Tutte's 2-factor Theorem\,(Theorem~\ref{tutte's theorem}),
$G$ has a barrier.
Let $(S, T)$ be a minimum barrier
such that $h_G(S,T)$ is smallest.
We use the same notations \textbf{$\CC_k$}, \textbf{$\CC_v$} and \textbf{$\CC_{1v}$} as defined in
the previous section.

\begin{CLA}\label{tutte}
$|T| > |S|+\sum_{k\ge 1}k\cdot |\CC_{2k+1}|$.
\end{CLA}
\spf
Since $(S, T)$ is a barrier,
\[
\begin{split}
\delta_G(S, T)&=2|S|-2|T|+\sum_{y\in T}d_{G-S}(y)-h_G(S, T)\\
&=2|S|-2|T|+\sum_{y\in T}d_{G-S}(y)-\sum_{k\ge 0} |\CC_{2k+1}| < 0.
\end{split}
\]
Let $U=V(G)-(S\cup T)$, by Lemma~\ref{minimal_barrier}~(1) and~(2),
\[
\sum_{y\in T}d_{G-S}(y)=\sum_{y\in T}e_G(y, U)
=e_G(T, U)=\sum_{k\ge 0}(2k+1)|\CC_{2k+1}|.
\]
Therefore,
we have
\[
0 > 2|S|-2|T|+\sum_{k\ge 0}(2k+1)|\CC_{2k+1}|-\sum_{k\ge 0}|\CC_{2k+1}|,
\]
which yields
$|T| > |S|+\sum_{k\ge 1}k|\CC_{2k+1}|$.
\qed

Based on the minimum barrier $(S,T)$, we define two bipartite graphs
$H^*$ and $H$ associated with $(S,T)$ as follows.  The definitions of $H^*$ and $H$ are fixed
hereafter.  The bipartite graph $H^*$ is defined as:
$$
V(H^*)=X'\cup T \quad \mbox{where}\quad X':=V(G)-T, \quad \mbox{and} \quad E(H^*)=E_G(X',T).
$$
Notice that the vertices from the even components in $G-(S\cup T)$\,(if any) are isolated
vertices in $H^*$ by  (2) of Lemma~\ref{minimal_barrier}.
The bipartite graph $H$ is obtained by performing the following operations to $G$.
\begin{enumerate}
\item
Remove all even components and all odd components in $\CC_1$.
\item
Remove all edges in $G[S]$.
\item
For a component $C\in\CC_{2k+1}$ with $k\ge 1$,
contract $C$ into one vertex and then split the
resulted vertex into
$k$ independent vertices
$U^C=\{u_1^C, u_2^C,\dots, u_k^C\}$
such that $d_H(u_1^C)=3$, and $d_H(u_2^C)=d_H(u_3^C)=\cdots=d_H(u^C_{k})=2$.
We note that,  this operation (3)
does  nothing to each single vertex component
$C\in \CC_3$.
\end{enumerate}

Let
$$
U^{\CC}=\bigcup_{k\ge 1}\left(\bigcup_{C\in\CC_{2k+1}}U^C\right)\quad \mbox{and}\quad X:=S\cup U^{\CC}.
$$
By the constructions,
the bipartite graph $H$ satisfies the following properties.

\begin{enumerate}
\item
$H$ is a bipartite graph with partite sets
$X$ and $T$,
\item
$|X|=|S|+\sum_{k\ge 1} k|\CC_{2k+1}|$,
and
\item
For each $k\ge 1$ and each $C\in\CC_{2k+1}$,
$d_H(u_1^C)=3$ and $d_H(u_i^C)=2$
for each $i$ with $2\le i\le k$.
\end{enumerate}
\par

Note that the construction of $H$ here is a modification of
the bipartite graph constructed in~\cite{local-CE}.
We now introduce some additional  notations. Those notations are used heavily in the subsequent proofs.

For each nonnegative integer $t$, let
$$
\CC_{\ge (2t+1)}:=\bigcup_{k\ge t} \CC_{2k+1}.
$$
It is clear that
$\CC_{\ge (2t+1)}\subseteq \HH_G(S,T)$.
For each $\mathcal{D}\subseteq \HH_G(S,T)$,  let
$$
V(\mathcal{D}):=\cup_{C\in \mathcal{D}} V(C),\quad
\mathcal{D}^1=\{C\in\mathcal{D}\,|\, |V(C)|=1\}\quad \mbox{and} \quad \mathcal{D}^2=\mathcal{D}-\mathcal{D}^1, \quad \mbox{and}\quad
U^{\mathcal{D}}=\left(\bigcup_{C\in\mathcal{D}}U^C\right).
$$
For example, we can take $\mathcal{D}=\CC_{\ge 3}\subseteq \HH_G(S,T)$  in the  above definition.
Then $V(\CC_{\ge 3})$ is the vertex set of all components $C\in\HH_G(S,T)$ such that
$e_G(T,V(C))\ge 3$; $\CC_{\ge 3}^1$ is the collection of  components
$C\in \CC_{\ge 3}$ such that
$|V(C)|=1$; $\CC_{\ge 3}^2$ is the collection of  components
$C\in\CC_{\ge 3}$ such that
$|V(C)|\ge 2$; and  $U^{\CC_{\ge 3}}$ is the set of vertices resulted
by splitting each contracted component in $\CC_{2k+1}$ into $k$ vertices,
for each integer $k\ge 1$.

Denote
$$
S':=S\cup V(\CC^1_{3}).
$$

\begin{CLA}\label{fact}
Each of the following holds:
\begin{itemize}
\item[$(1)$] $d_{H^*}(y)=d_G(y)$ for each $y\in T$;
\item[$(2)$] $X'\cap X=S'=S\cup V(\CC_3^1)$;
\item[$(3)$] $d_H(x)=d_{H^*}(x)$ for each $x\in S'$;
 \item[$(4)$] $d_H(y)=d_{H^*}(y)-|\CC_{1y}|$  for each $y\in T$.
\end{itemize}
\end{CLA}

\spf The  statements (1)-(3) are obvious. We only show  the last one.
By (2) of Lemma~\ref{minimal_barrier}
that for each $y\in T$
and each even component $C$ of
$G-(S\cup T)$, $e_G(y,V(C))=0$ holds. Thus
 $d_G(y)=e_G(y,S)+e_G(y,V(\CC_{\ge 1}))=
 e_G(y,S\cup V(\CC_{\ge 3}))+e_G(y,V(\CC_{1}))=d_{H^*}(y)$.
By the construction of
$H$, $d_{H}(y)=e_G(y,S\cup V(\CC_{\ge 3}))=d_G(y)-e(y,V(\CC_1))=d_{H^*}(y)-|\CC_{1y}|$.
\qed

In the remaining proof, using Lemma~\ref{matching},
we first show that, in the bipartite graph
$H^*$,  there is a matching which saturates
$T$. Then by using the relations between $H^*$ and
$H$ and Hall's Theorem, we show that, in $H$,
there is a matching  which saturates
$T$. The later one gives that $|X|=|S|+\sum_{k\ge1}k|\CC_{2k+1}|\ge |T|$,
leading  a contradiction to
Claim~\ref{tutte}.



\begin{CLA}\label{noDelta}
$T$ has no $\Delta$-vertex.
\end{CLA}

\spf Suppose on the contrary that there
exists $w\in T$ such that $d_G(w)=\Delta$.
Let $V_{\text{even}}$ be the  vertex set of
the even components in $G-(S\cup T)$.
Then
 by  $e_G(w,V(\CC_{\ge1}))= |\CC_{1w}|+|\CC_{w}-\CC_{1w}|$\,((3) of Lemma~\ref{minimal_barrier}) and
$|X|=|S|+\sum_{k\ge 1}k|\CC_{2k+1}|<|T|$\,(Claim~\ref{tutte}), 
\begin{eqnarray*}
 &&\frac{1}{2}\left(|S|+|T|+|V(\CC_{\ge 1})|+|V_{\text{even}}|\right)= \frac{n}{2}\le  \Delta=d_G(w) \le  |S|+|\CC_{1w}|+|\CC_{w}-\CC_{1w}| \\
 &\le & \frac{1}{2}(|S|+2|\CC_{1w}|+|\CC_{w}-\CC_{1w}|+(|S|+\sum\limits_{k\ge 1}|\CC_{2k+1}|))\quad \mbox{(by $\CC_w-\CC_{1w}\subseteq \bigcup_{k\ge 1}\CC_{2k+1}$)}\\
   &< & \frac{1}{2}\left(|S|+2|\CC_{1w}|+|\CC_{w}-\CC_{1w}|+|T|\right).
\end{eqnarray*}
The above strict inequalities give that
\begin{eqnarray}
2|\CC_{1w}|+|\CC_{w}-\CC_{1w}| &\ge  & |V(\CC_{\ge 1})|+|V_{\text{even}}|+1 \nonumber \\
  & \ge &(|V(\CC_{1w})|+1) +\sum_{C\in \CC_1-\CC_{1w}} |V(C)|+|V(\CC_{\ge 3})|+|V_{\text{even}}|\label{comp}.
  \end{eqnarray}
Since  $|V(\CC_{\ge 3})|\ge |\CC_{w}-\CC_{1w}|$, we have that $2|\CC_{1w}|\ge (|V(\CC_{1w})|+1) +\sum_{C\in \CC_1-\CC_{1w}} |V(C)| +|V_{\text{even}}|$.
If $|\CC_{1w}|\ge 2$, then $|V(\CC_{1w})|\ge 2|\CC_{1w}|$\,(Lemma~\ref{size2}), showing a contradiction.
If $|\CC_{1w}|=0$, then $|V(\CC_{1w})|=0$. So
$|\CC_{1w}|=1$.  Since $|V(\CC_{1w})|\ge |\CC_{1w}|=1$ and $|V(\CC_{1w})|+1\le 2$,
we get that $|V(\CC_{1w})|=1$ and so \textbf{(a)}\,$\CC_1=\CC_{1w}$ and  $|V_{\text{even}}|=0$.
Using $|V(\CC_{\ge 3})|\ge |\CC_{w}-\CC_{1w}|$ again, under the above facts,
inequality~(\ref{comp}) becomes \textbf{(b)}\,$ |\CC_{w}-\CC_{1w}| =|V(\CC_{\ge 3})|$.

%
Then (a),
 together with the fact  that $|V(\CC_{1w})|=1$ implies that there is exact one single vertex component in $\CC_1$.
By $|V_{\text{even}}|=0$ in (a), there is no even component in $G-(S\cup T)$.
Since $|V(\CC_{1w})|=1$, by Lemma~\ref{size2}, $\CC_{w}=\CC_{1w}$.
By $|\CC_{w}-\CC_{1w}|=|V(\CC_{\ge 3})|$ in (b),  we see that
$V(\CC_{\ge 3})=\emptyset$.  Let $H^*[X',T]$ and $H[X,T]$ be the two
bipartite graphs associated with $(S,T)$.  Then  $|X'-X|=|V(\CC_1)|=|V(\CC_{1w})|=1$.
Combining $|X'|\ge d_{H^*}(w)=d_G(w)=\Delta\ge n/2$,
we have $|X|\ge n/2-1$. As $|T|>|X|$\,(Claim~\ref{tutte}) and $|T|+|X'|=|T|+|X|+1=n$,
we get $|T|=|X'|=n/2$ and $|X|=n/2-1\le \Delta-1$.
Hence, for any $y\in T\cap V_{\Delta}$, $y$ is adjacent to the
unique  vertex in $\CC_{1w}=\CC_1$.
As $e_G(T,V(\CC_{1w}))=1$,  $T\cap V_{\Delta}=\{w\}$.
That is, $w$ is the unique  $\Delta$-vertex in $T$.
Applying
Lemma~\ref{TnoDelta} with $\delta_0=1$ on $H^*$,
we see that for every
edge
$xy\in E( H^*)$ with $x\in X'$ and $y\in T$, the
relation  $d_{H^*}(y)\ge d_{H^*}(x)$ holds.
As $d_{H^*}(y)=d_G(y)\ge 2$,
 $T$ has no isolated vertices.
Hence in $H^*$, there is a matching $M$ which saturates $T$.
Since $|X'|=|T|$,
 $M$ is a perfect matching of $H^*$.
%
 Let
$w^*$ be the vertex to which $w$ is
adjacent in $\CC_1$. Then $V(M)$ contains $w^*$
and $d_{H^*}(w^*)=1$ by noticing that $d_{H^*}(w^*)=e_G(w^*,T)=e_G(T, V(\CC_{1w}))=1$.
As for any $y\in T$,
$d_{H^*}(y)=d_G(y)\ge 2$,
$d_{H^*}(w)>d_{H^*}(w^*)$.
For any $xy\in E( H^*)-\{ww^*\}$ with $x\in X'$ and $y\in T$,
 $d_{H^*}(y)\ge d_{H^*}(x)$.
Hence,
\begin{eqnarray*}
  e_{H^*}(X', T) &=  & \sum\limits_{xy\in M-\{ww^*\}\atop x\in X', y\in T}d_{H^*}(y)+d_{H^*}(w)  \\
   & > & \sum\limits_{xy\in M-\{ww^*\}\atop x\in X', y\in T}d_{H^*}(x)+d_{H^*}(w^*)=  e_{H^*}(X', T),
\end{eqnarray*}
showing  a contradiction.
%
\qed

For a vertex $x\in X'=V(G)-T$, define $\sigma_x$
as the number of non $\Delta$-degree neighbors of
$x$ in $X'$ and let
$$
S_1=\{x\in S'=S\cup V(\CC_{3}^1)\,|\,\sigma_x\ge 1\} \quad \text{and} \quad
S_0=S'-S_1=\{x\in S'\,|\sigma_x=0\}.
$$
Following the definitions of $S_0$ and $S_1$,  $N_G(x)\cap S_0=\emptyset$ for any non $\Delta$-degree vertex $x\in X'$.


As $T$ has no $\Delta$-vertex, applying
Lemma~\ref{TnoDelta} with $\delta_0=0$, for each edge $xy\in E(H^*)$
with $x\in X'$ and $y\in T$, we have the following claim.

\begin{CLA}\label{Hstardegree}
 $$d_{H^*}(y) \ge
\left\{
  \begin{array}{ll}
    d_{H^*}(x)+2, & \hbox{if $x\in S_1$;} \\
    d_{H^*}(x)+1, & \hbox{otherwise.}
  \end{array}
\right.
$$
Moreover, $H^*$ has a matching which saturates $T$.
\end{CLA}

\begin{CLA}\label{Delta-largestorder}
Let $C\in \CC_{\ge1}$ and  $x\in V(C)\cap V_{\Delta}$.
If $e_G(x,T)\le 1$,
then $|V(C)|>\frac{1}{2}|V(\CC_{\ge 1})|$.
\end{CLA}

\spf Suppose on the contrary that
$|V(C)|\le \frac{1}{2}|V(\CC_{\ge 1})|$, that is,
$|V(C)|\le |V(\CC_{\ge1})|-|V(C)|$.
Then since $e_G(x,T)\le1$,
 $$n/2\le \Delta= d_G(x)\le |S|+|V(C)|-1+e_G(x,T)\le |S|+|V(C)|\le |S|+|V(\CC_{\ge1})|-|V(C)|.
 $$
 As $|T|>|S|$,
 $$
 n\le 2|S|+|V(C)|+|V(\CC_{\ge1})|-|V(C)|<|S|+|T|+|V(C)|+|V(\CC_{\ge1})|-|V(C)|\le n,
 $$
showing a contradiction.
\qed

For each $y\in T$, we have $d_{H^*}(y)=d_G(y)\ge 2$. The following claim gives
a property when $d_{H^*}(y)=2$.

\begin{CLA}\label{degreeHstarC3}
Let  $y\in T$ be a vertex. If 
$d_{H^*}(y)=2$, then there exists $x\in N_{H^*}(y)\cap S$
such that $d_{H^*}(x)=1$.
\end{CLA}

\spf
By Vizing's Adjacency Lemma,
each vertex of $G$ is adjacent to
at least two $\Delta$-vertices.
 Since $d_{H^*}(y)=d_G(y)=2$,
 the two neighbors of $y$
are $\Delta$-vertices.
By Claim~\ref{Hstardegree},
$2=d_{H^*}(y)\ge d_{H^*}(x)+1$, so
 each of the
two neighbors of $y$ has degree 1 in $H^*$\,(
and thus has exact one neighbor in $T$).
If $N_{H^*}(y)\subseteq V(\CC_{\ge1})$,
then by (3) of Lemma~\ref{minimal_barrier},
each of the vertex in $N_{H^*}(y)$
is contained in a distinct component in
$\CC_{\ge 1}$.
However, by Claim~\ref{Delta-largestorder},
there exists at most one component
$C\in \CC_{\ge 1}$ such that
it contains a $\Delta$-vertex and the
$\Delta$-vertex has degree exact 1 in $H^*$.
Hence $N_{H^*}(y)\cap S\ne \emptyset$.
Let $x\in N_{H^*}(y)\cap S$. Then
$x$ is the desired vertex.
\qed



\begin{CLA}\label{c1_ne_empty}
We may assume that $\CC_{1}\ne \emptyset$.
\end{CLA}

\spf Suppose on the
contrary that $\CC_{1}=\emptyset$.
By Claim~\ref{fact},
for each $x\in S$, $d_H(x)=d_{H^*}(x)$ and
 for any
$y\in T$, $d_H(y)=d_{H^*}(y)$.
Applying  Claim~\ref{degreeHstarC3},
if
there exists $y\in T$ such that $d_H(y)=2$, then $y$
has a neighbor of degree 1 in $S$.
Let $X_1=\{x\in X\,|\,d_H(x)=1\}$ and $T_1=N_H(X_1)$,
and let $H'[X-X_1, T-T_1]:=H[(T-T_1)\cup (X-X_1)]$.
Then for any $y\in T-T_1$, $d_{H'}(y)\ge 3$.
We then claim that for each edge $xy\in E(H')$
with $x\in X-X_1$ and $y\in T-T_1$, $d_{H'}(y)\ge d_{H'}(x)$
holds.
If $x\in S$, then by Claim~\ref{Hstardegree}, $d_{H^*}(x)\le d_{H^*}(y)$.
Hence $d_{H'}(x)\le d_H(x)=d_{H^*}(x)\le d_{H^*}(y)=d_H(y)=d_{H'}(y)$.
If $x\in U^{\CC}$, then $d_{H'}(x)\le d_H(x)\le 3$ by the construction of $H$.
Hence $d_{H'}(x)\le 3\le d_{H'}(y)$.
Since $N_H(X_1)=T_1$,  that $T$ has no isolated vertices in $H$ implies that
$T-T_1$ has no isolated vertices in $H'$.
By Lemma~\ref{matching},
$H'$ has a matching which saturates $T-T_1$.
By Lemma~\ref{islatedmatching}, $H$ has a matching which
 which saturates $T$. This
 gives a contradiction to Claim~\ref{tutte}.
\qed

%


For each $C\in \CC_1$, by (3) of Lemma~\ref{minimal_barrier},
there is a unique vertex $y_c\in T$  adjacent to
a unique vertex $x_c$ on it.
We call $x_c$ and $y_c$ the partners of each other.
We divide the components in $\CC_1$ into two subgroups $C_{11}$ and $C_{12}$
in order to consider the degrees of the
partner vertices in $T$:
$$
\CC_{11}=\{C\in \CC_1\,|\, e_G(y_{c}, \CC_1)=1\}\quad \text{and }\quad \CC_{12}=\{C\in \CC_1\,|\, e_G(y_{c}, \CC_1)\ge 2\}.
$$
By the definition of $\CC_{12}$, it is clear that if $\CC_{12}\ne \emptyset$,
then $|\CC_{12}|\ge 2$.
Also, by Lemma~\ref{size2}, for each $C\in \CC_{12}$, $|V(C)|\ge 2$.

Furthermore, we divide the  components in $\CC_{11}$ into
two groups as follows.
$$
\CC_{11}^1=\{C\in \CC_{11}\,\mid\, |V(C)|=1\}\quad \text{and} \quad \CC_{11}^2=\{C\in \CC_{11}\,\mid\, |V(C)|\ge 2\}.
$$
Corresponding to the partition of
$\CC_1$, we partition vertices in $T$
into subgroups, as follows.
\begin{eqnarray*}
  T_{1}^1&=& \{y\in T\,|\, e_G(y,V(\CC_{11}^1))=1\},\quad T_{1}^2=\{y\in T\,|\, e_G(y,V(\CC_{11}^2))=1\}; \\
  T_0&=& \{y\in T\,|\, e_G(y,V(\CC_{1}))=0\},\quad \text{and}\quad T_2=\{y\in T\,|\, e_G(y,V(\CC_{12}))\ge 2\}.
\end{eqnarray*}
Notice that for a vertex $y\in T$, if $e_G(y,V(\CC_{1}))\ge2$, then $e_G(y,V(\CC_{1}))=e_G(y,V(\CC_{12}))$.
Hence,
\begin{equation}\label{T12}
 e_G(y, V(\CC_1))=1 \quad  \mbox{for each}\quad y\in T_1^1\cup T_1^2.
\end{equation}
Since each $C\in \CC_{11}^1$ satisfies $|V(C)|=1$, by Lemma~\ref{size2},
\begin{equation}\label{T1}
 e_G(y, V(\CC_{\ge 1}))=1 \quad  \mbox{for each}\quad y\in T_1^1.
\end{equation}

Let \quad\quad  $m_{11}:=|\CC_{11}^1|, \quad m_{12}:=|\CC_{11}^2|,\quad  m_{2}:=|\CC_{12}|, \quad \text{and} \quad m_3:=|\CC^2_{\ge3}|.$


\begin{CLA}\label{c111_no_Delta}
We may assume that none vertices in $V(\CC_{11}^1)$
is a $\Delta$-vertex.
\end{CLA}

\spf Suppose on the country and let  $x_c\in V(\CC_{11}^1)\cap V_{\Delta}$.
Since $e_G(x_c,T)=1$ and $e_G(x_c,V(G)-S-T)=0$,
we have $e_G(x_c,S)=\Delta-1\ge n/2-1$.
This indicates that $|S|\ge n/2-1$.
Combining $|T|>|S|+\sum_{k\ge1}k|\CC_{2k+1}|$\,(Claim~\ref{tutte})
and $|S|+|T|< |X'|+|T|=n$\,(noticing that $|S|<|X'|$ by $1=|\{x_c\}|\le |V(\CC_{11}^1)|$ and $V(\CC_{11}^1)\cup S\subseteq X'$), we have
$|T|=n/2=|S|+1$.   We consider the
bipartite graph $H^*[X',T]$ associated with $(S,T)$.
As $|V(G)|=n$ and $|T|=n/2$,
$|X'|=|T|$.
By Claim~\ref{Hstardegree}, $H^*$
has a matching $M$ which saturates $T$.
Since $|T|=|X'|$, $M$ is a
perfect matching.
Since  $d_{H^*}(x_c)=1$,
$T$ has a unique neighbor, say $y_c$ of
$x_c$.
Then $x_cy_c\in M$.
Because $d_{H^*}(y_c)=d_G(y_c)\ge 2$,
$d_{H^*}(y_c)>d_{H^*}(x_c)$.
By Claim~\ref{Hstardegree},
$d_{H^*}(y)\ge d_{H^*}(x)$
for each $xy\in E(H^*)-\{y_cx_c\}$ with $x\in X'$ and
$y\in T$. Hence,
\begin{eqnarray*}
  e_{H^*}(X', T) &=  & \sum\limits_{xy\in M-\{x_cy_c\}\atop x\in X', y\in T}d_{H^*}(y)+d_{H^*}(y_c)  \\
   & > & \sum\limits_{xy\in M-\{x_cy_c\}\atop x\in X', y\in T}d_{H^*}(x)+d_{H^*}(x^c) = e_{H^*}(X', T),
\end{eqnarray*}
showing a contradiction.
\qed

%

By the definition, for each $ C\in \CC^2_{\ge 3}\cup \CC_{11}^2\cup \CC_{12}$,
we have $|V(C)|\ge 2$ holds.  Thus
\begin{eqnarray}\label{size_n}
   n &\ge & |S'|+|T|+|V(\CC_{11}^1)|+|V(\CC_{11}^2)|+|V(\CC_{12})| +|V(\CC^2_{\ge 3})|\\
   &\ge &  |S'|+|T|+m_{11}+2m_{12}+2m_{2}+2m_3,\nonumber
\end{eqnarray}
where   $S'=S\cup V(\CC_{ 3}^1)$ is defined previously.

\begin{CLA}\label{oddatleast5}
Let $xy\in E(H^*)$ be an edge with
$x\in V(C)\subseteq V(\CC_{\ge 5}^1)$ and $y\in T$,
and  let $u_c$ be a vertex in $U^{C}$
which is adjacent to $y$ in $H$.
Then  $d_{H^*}(y)\ge  d_H(u_c)+3 $.
\end{CLA}

 \spf  Let
$V(C)=\{x\}$. Then $d_{H^*}(x)\ge 5$ as $C\in \CC_{\ge 5}$.
  By Claim~\ref{Hstardegree}, $d_{H^*}(y)\ge 6$.
  Recall that in $U^C$, $d_H(u^C_1)=3$ and
  $d_H(u^C_i)=2$ for $i\ge 2$, so
  $d_{H^*}(y)\ge 6\ge d_H(u_c)+3 $.
\qed

For each vertex $y\in T_1^1\cup T_1^2\cup T_2$,
$|\CC_{1y}|\ge 1$. So $d_H(y)=d_{H^*}(y)-|\CC_{1y}|<d_{H^*}(y)=d_G(y)$.
In order to find a matching saturating $T$ in $H$,
in the following three claims, we show that $y$
still has enough neighbors remained in $V(G)-T-V(\CC_1)=S\cup V(\CC_{\ge 3})$.

\begin{CLA}\label{c111}
If $T_1^1\ne \emptyset$, then for each
$y\in T_1^1$,
$$
|N_G(y)\cap V_{\Delta}\cap S|\ge|S_0|+(m_{11}+1)/2+m_{12}+m_2+m_3.
$$
\end{CLA}

\spf Let $x\in V(\CC_{11}^1)$ such that $xy\in E(G)$.
Since $y\in T_1^1$, using (\ref{T1}) that $e(y,V(\CC_{\ge 1}))=1=e(y,V(\CC_{11}^1))$,
$N_G(y)\cap V(\CC_{\ge 1})=\{x\}$.
By Claim~\ref{c111_no_Delta},
$x$ is not a $\Delta$-vertex.
Thus $N_G(y)\cap V_{\Delta}\subseteq S$. So we only need to show
$|N_G(y)\cap V_{\Delta}|\ge|S_0|+(m_{11}+1)/2+m_{12}+m_2+m_3$.
Recall $S_0=\{x\in S'\,|\,\sigma_x=0\}$
is the set of vertices in $S'=S_0\cup S_1$ only adjacent to $\Delta$-vertices
in $X'=V(G)-T$, so $N_G(x)\cap S_0=\emptyset$.
Hence $d_G(x)\le |N_G(x)\cap T|+|N_G(x)\cap S_1|\le |S_1|+1$.
By Vizing's Adjacency Lemma, $y$ is adjacent to
at least $\Delta-d_G(x)+1$
$\Delta$-vertices in $G$.
Simple calculation shows that
\begin{eqnarray*}
  \Delta-d_G(x)+1 &\ge  & n/2- |S_1|-1+1\\
   &\ge  & \frac{1}{2}(|S' |+|T|+m_{11}+2m_{12}+2m_2+2m_3)-|S_1|\quad \mbox{(by inequality~(\ref{size_n}))} \\
   &\ge  &|S_0|+(m_{11}+1)/2+m_{12}+m_2+m_3,
\end{eqnarray*}
where the last inequality is obtained by using the facts that $|S'|=|S_1|+|S_0|$ and $|S'|
=|S|+|V(\CC_3^1)|=|S|+|\CC_3^1|\le |S|+\sum_{k\ge1}|\CC_{2k+1}|<|T|$ by Claim\ref{tutte}.
\qed

If $\CC_{11}^2\ne \emptyset$,
let $C_{\text{max}}^1\in \CC_{11}^2$ be a component
such that $|V(C_{\text{max}}^1)|=\max\{|V(C)|\,|\,C\in \CC_{11}^2\}$.
Then by Claim~\ref{Delta-largestorder}, if $V(\CC_{11})\cap V_{\Delta}\ne \emptyset$,
then $V(\CC_{11})\cap V_{\Delta}=V(C_{\text{max}}^1)\cap V_{\Delta}$\,
(since $\CC_{11}=\CC_{11}^1\cup \CC_{11}^2$ and $V(\CC_{11}^1)\cap V_{\Delta}=\emptyset$
by Claim~\ref{c111_no_Delta}).

\begin{CLA}\label{c112}
If $T_1^2\ne \emptyset$, then for each
$y\in T_1^2$,
$$|N_G(y)\cap V_{\Delta}\cap (S\cup V(\CC_{\ge 3}))|\ge
\left\{
  \begin{array}{ll}
    |S_0|+(m_{11}+1)/2+m_{12}+m_2+m_3-1, & \hbox{if  $x\notin V(C_{\text{max}}^1)$;} \\
    1, & \hbox{if $x\in V(C_{\text{max}}^1)$;}
  \end{array}
\right.
$$
where $x$ is the neighbor
of $y$ in $\CC_{11}^2$.
\end{CLA}

\spf  
Suppose first that $C\neq C_{\text{max}}^1$.
Since $y$ is adjacent to exactly one
component $C$ in $\CC_{11}^2$\,(by~(\ref{T12})) and
$C$ contains no $\Delta$-vertex,  $N_G(y)\cap V_{\Delta}\subseteq S\cup V(\CC_{\ge 3})$.
So we only need to show
$|N_G(y)\cap V_{\Delta}|\ge|S_0|+(m_{11}+1)/2+m_{12}+m_2+m_3-1$.
Again, as $x\notin V_{\Delta}$,  $N_G(x)\cap S_0=\emptyset$.
Hence $d_G(x)\le |N_G(x)\cap T|+|N_G(x)\cap S_1|+|V(C)|-1\le |S_1|+|V(C)|$.
By Vizing's Adjacency Lemma, $y$ is adjacent to
at least $\Delta-d_G(x)+1$
$\Delta$-vertices in $G$.
Since $n\ge |S'|+|T|+m_{11}+|V(C)|+|V(C_{\text{max}}^1)|+2(m_{12}-2)+2m_2+2m_3$,
\begin{eqnarray*}
  \Delta-d_G(x)+1 &\ge  & n/2-(|S_1|+|V(C)|)+1\\
   &\ge  & \frac{1}{2}(|S'|+|T|+m_{11}+|V(C)|+|V(C_{\text{max}}^1)|+2(m_{12}-2)+2m_2+2m_3)-|S_1|-|V(C)|+1 \\
   &\ge  &|S_0|+(m_{11}+1)/2+m_{12}+m_2+m_3-1.
\end{eqnarray*}

Suppose now that $C=C_{\text{max}}^1$.
As $d_G(y)\ge 2$, and $e_G(y,V(\CC_1))=1$ by~(\ref{T12}),
the other neighbor of
$y$ is contained in $S\cup V(\CC_{\ge 3})$.
\qed

If $C_{\text{max}}^1$ exists, let $y_s$\,(for \textbf{$y_{\text{special}}$}) be the unique
vertex in $T$ such that $e_G(y_s,V(C_{\text{max}}^1))=1$.

\begin{CLA}\label{c12}
If $T_2\ne \emptyset$, then for each
$y\in T_2$,
$$
|N_G(y)\cap V_{\Delta}\cap (S\cup V(\CC_{\ge 3}))|\ge|S_0|+m_{11}/2+m_{12}+m_2+m_3-1.
$$
\end{CLA}

\spf If $T_2\ne \emptyset$, then by the definition, $\CC_{12}\ne \emptyset$,
giving that $m_{2}\ge 2$.
Let $C^2_{\text{max}}$ be a component with largest cardinality in $\CC_{12}$.
Let  $C_1\in \CC_{12}\cap \CC_{1y}-\{C^2_{\text{max}}\}$
and $x$ be the neighbor of $y$ on $C_1$.
By Claim~\ref{Delta-largestorder},  if $C^2_{\text{max}}$ contains a $\Delta$-vertex,
it is the only
component in $\CC_{1}$  which contains a $\Delta$-vertex.
Thus, $N_G(y)\cap V_{\Delta}\subseteq S\cup V(\CC_{\ge 3})\cup V(C^2_{\text{max}})$.
So it suffices to show that
$|N_G(y)\cap V_{\Delta}|-|N_G(y)\cap V(C^2_{\text{max}})\cap V_{\Delta}|\ge|S_0|+m_{11}/2+m_{12}+m_2+m_3-1$.
Again, as $x\notin V_{\Delta}$,  $N_G(x)\cap S_0=\emptyset$.
Hence $d_G(x)\le |N_G(x)\cap T|+|N_G(x)\cap S_1|+|V(C)|-1\le |S_1|+|V(C)|$.
By Vizing's Adjacency Lemma, $y$
is adjacent to at least $\Delta-d_G(x)+1$ $\Delta$-vertices in $G$.
So $y$ has at least $\Delta-d_G(x)+1$ $\Delta$-degree
neighbors in $S\cup V(\CC_{\ge 3})$ if
$|N_G(y)\cap V(C^2_{\text{max}})\cap V_{\Delta}|=0$; and
$y$ has at least $\Delta-d_G(x)$ $\Delta$-degree
neighbors in $S\cup V(\CC_{\ge 3})$ if $|N_G(y)\cap V(C^2_{\text{max}})\cap V_{\Delta}|=1$.

If $|N_G(y)\cap V(C^2_{\text{max}})\cap V_{\Delta}|=0$,  we get that
\begin{eqnarray*}
  \Delta-d_G(x)+1 &\ge  & n/2-(|S_1|+|V(C_1)|)+1\\
   &\ge &   \frac{1}{2}(|S'|+|T|+m_{11}+2m_{12}+|V(C_1)|+|V(C^2_{\text{max})}|+2(m_2-2)+2m_3)-|S_1|-|V(C_1)|+1 \\
   &\ge  &|S_0|+(m_{11}+1)/2+m_{12}+m_2+m_3-1.
\end{eqnarray*}

If $|N_G(y)\cap V(C^2_{\text{max}})\cap V_{\Delta}|=1$,
then $C^2_{\text{max}}$ contains a $\Delta$-vertex $x$ with $e_G(x,T)=1$,
and thus $|V(C^2_{\text{max}})|>|V(C_1)|$ by Claim~\ref{Delta-largestorder}.
Also since $|S'|<|T|$, we get
\begin{eqnarray*}
  \Delta-d_G(x) &\ge  & n/2-(|S_1|+|V(C_1)|)\\
   &\ge  & \frac{1}{2}(|S'|+|T|+m_{11}+2m_{12}+|V(C_1)|+|V(C^2_{\text{max}})|+2(m_2-2)+2m_3)-|S_1|-|V(C_1)| \\
   &\ge  &|S_0|+m_{11}/2+m_{12}+m_2+m_3-1.
\end{eqnarray*}
\qed

\begin{CLA}\label{T-no-isolated}
In the bipartite graph $H[X,T]$,
$T$ has no isolated vertices.
\end{CLA}

\spf Let $y\in T$ be a vertex.
If $|\CC_1|\le 1$, then $d_H(y)\ge d_{H^*}(y)-1=d_G(y)-1\ge 1$.
So assume $|\CC_1|\ge 2$.
If $y\in T_0$,
then $d_H(y)=d_{H^*}(y)=d_G(y)\ge 2$.
For each $y\in T_1^1\cup T_1^2\cup T_2$,
either $d_H(y)\ge m_{11}/2+m_{12}+m_2-1$ or $d_H(y)\ge 1$
by claims~\ref{c111}-\ref{c12}.
Since $m_{11}+m_{12}+m_2=|\CC_1|\ge 2$,
$d_H(y)\ge 1$.
\qed.

\begin{CLA}\label{Hdegree}
Let $xy\in E(H)$ be an edge with $x\in X$ and $y\in T$.
Then each of the following
holds:
\begin{itemize}
  \item [$(1)$] $d_H(y)+|\CC_{1y}|\ge d_H(x)+2$ if $x\in S_1$;
  \item [$(2)$] $d_H(y)+|\CC_{1y}|\ge d_H(x)+1$ if $x\in S_0$;
   \item [$(3)$]$d_H(y)+|\CC_{1y}|\ge d_H(x)+3$ if $x\in U^{\CC_{\ge 5}^1}$;
   \item [$(4)$]For each $x\in U^{\CC_{\ge 3}^2}$, either
   $d_H(y)+|\CC_{1y}|\ge d_H(x)$ or $d_H(y)+|\CC_{1y}|=2$ and $d_H(x)=3$.
   In the later case, there exists $x\in S$ such that $xy\in E(H)$ and $d_H(x)=1$.
\end{itemize}
\end{CLA}

\spf As $d_H(y)+|\CC_{1y}|=d_{H^*}(y)$ and $d_H(x)=d_{H^*}(x)$  for all $x\in S'=X\cap X'$,
(1)-(2) follow Claim~\ref{Hstardegree}.
By Claim~\ref{oddatleast5}, we get (3).
For each $y\in T$, $d_{H^*}(y)=d_H(y)+|\CC_{1y}|=d_G(y)\ge 2$, and for each
$x\in U^{\CC_{\ge 3}^2}$, according to the construction of $H$, either $d_H(x)=2$
or $d_H(x)=3$. If $d_{H^*}(y)\ge 3$, the first part of (4) holds.
If $d_{H^*}(y)=2$, then the second part of (4) follows. The existence of the vertex
$x\in S$ such that $xy\in E(H)$ and $d_H(x)=1$ is guaranteed by
Claim~\ref{degreeHstarC3}.
\qed

Let $y\in T$ be a vertex of degree 2 in $H^*$.
By Claim~\ref{degreeHstarC3}, $y$
has a neighbor $x$ in $S$ which has degree
1 in $H^*$. As $S\subseteq X\cap X'$, $y$ has a neighbor $x$
of degree 1 also in $H$.
Applying Lemma~\ref{islatedmatching},
to show that $H$ has a matching which saturates $T$,
we may assume that
for any vertex $y\in T$,
$d_{H^*}(y)\ge 3$ holds.
By Claim~\ref{Hdegree},  the assumption indicates that
\begin{equation}\label{degreeHrefine}
  d_H(y)+|C_{1y}|\ge 3 \quad \mbox{and}\quad d_H(y)+|\CC_{1y}|\ge d_H(x) \,\, \mbox{for every edge }\,\, xy\in E(H).
\end{equation}


\begin{CLA}\label{Hmatching}
$H$ has a matching which saturates $T$.
\end{CLA}

\spf
Suppose not. Then by Hall's Theorem, there
is a nonempty set $A\subseteq T$ such that
$|N_H(A)|<|A|$. We choose  $A$ such that
it has the minimum cardinality and satisfies
$|N_H(A)|<|A|$. Let $B:=N_H(A)$ and $H':=H[A\cup B]$.
We claim that, in $H'$,  there is a matching which saturates
$B$. Suppose on the contrary. Then
by Hall's Theorem again, there is a nonempty subset $B'\subseteq B$
such that $|N_{H'}(B')|<|B'|$.
Since $B'\subseteq B=N_H(A)\ne \emptyset$\,($T$ has no isolated vertices by
Claim~\ref{T-no-isolated}), $N_{H'}(B')\ne \emptyset$.  Let $A'=A-N_{H'}(B')$.
As $|A|>|B|\ge |N_{H'}(B')|>0$, $0<|A'|<|A|$.
On the other hand, we have $N_{H'}(A')=N_H(A')=B-B'$.
However, $|A'|=|A|-|N_{H'}(B')|>|B|-|N_{H'}(B')|> |B|-|B'|=|B-B'|=|N_H(A')|$,
showing a contradiction to the choice of $A$.

In $H'$, let $M$ be a matching  which saturates $B$.
We consider three cases below.

{\noindent \emph{Case 1. }} $A \subseteq T_0$.

In this case, all vertices  $y\in A$
has $|\CC_{1y}|=0$.  By Claim~\ref{Hdegree},
 $d_{H'}(y)=d_H(y)\ge d_{H'}(x)$
for every edge $xy\in E(H')$.
As $|A|>|B|$ and $M$ saturates $B$, $A-V(M)\ne \emptyset$.
Let $y^*\in A-V(M)$.
Then $d_H(y^*)\ge 3$ by ~(\ref{degreeHrefine}).
Then we get
\begin{eqnarray*}
  e_H(A,B) &=& \sum\limits_{xy\in M\atop x\in B, y\in A} d_{H'}(x)\le \sum\limits_{xy\in M\atop x\in B, y\in A} d_{H}(x) \\
    &\le & \sum\limits_{xy\in M\atop x\in B, y\in  A} (d_H(y)+|\CC_{1y}|)\quad\mbox{(by~(\ref{degreeHrefine}))}\\
 & < &  \sum\limits_{xy\in M\atop x\in B, y\in  A} d_H(y)+d_H(y^*)\le e_H(A,B)\quad\mbox{($|\CC_{1y}|=0$ for $ y\in A$ and $d_{H}(y^*)\ge 3$)},
\end{eqnarray*}
showing a contradiction.

{\noindent \emph{Case 2. }} $A\cap (T_1^1\cup T_1^2\cup T_2)=A\cap T_1^2=\{y_s\}$.

Note that in this case, $y_s\in T_1^2$ and  thus $|\CC_{1y_s}|=e_G(y_s,V(\CC_1))=1$ by~(\ref{T12}).
Since $A\cap(T_1^1\cup T_1^2\cup T_2)=\{y_s\}$, for each $y\in A-\{y_s\}$,
$|\CC_{1y}|=0$.  As $|A|>|B|$ and $M$ saturates $B$, $A-V(M)\ne \emptyset$.
Let $y^*\in A-V(M)$.
Following~(\ref{degreeHrefine}), we have that $d_{H^*}(y)+|\CC_{1y^*}|\ge 3$.
If $y^*\ne y_s$, $d_H(y^*)=d_{H}(y^*)+|\CC_{1y^*}|\ge 3$;
and  if  $y^*= y_s$ then $|\CC_{1y^*}|=1$, so $d_H(y^*)\ge 2$. We may assume that $y_s\in A\cap V(M)$,
for otherwise, we can get a contradiction by the same argument as in Case 1. Hence,
\begin{eqnarray*}
 && e_H(A,B) = \sum\limits_{xy\in M\atop x\in B, y\in A} d_{H'}(x)\le \sum\limits_{xy\in M\atop x\in B, y\in A} d_{H}(x) \\
    &\le & \sum\limits_{xy\in M\atop x\in B, y\in A-\{y_s\}} (d_H(y)+|\CC_{1y}|)+(d_H(y_s)+1)\quad\mbox{(by~(\ref{degreeHrefine}))}\\
 & \le &  \sum\limits_{xy\in M\atop x\in B, y\in A} d_H(y)+1+(d_H(y^*)-3)\le e_H(A,B)-2\quad\mbox{($|\CC_{1y}|=0$ for $ y\in A-\{y_s\}$ and $d_H(y^*)\ge 3$)},
\end{eqnarray*}
giving a contradiction.

{\noindent \emph{Case 3. }} $A\cap(T_1^1\cup T_1^2\cup T_2)-\{y_s\}\ne \emptyset$.

Let $y'\in A\cap (T_1^1\cup T_1^2\cup T_2)-\{y_s\}$ such that
$|N_G(y')|=\max\limits_{y\in A\cap (T_1^1\cup T_1^2\cup T_2)-\{y_s\}}|N_G(y)|$.
Denote $B_1:=N_H(y')$ and 
 $\overline{B_1}:=B-B_1$.  Then $V(M)\cap B=B=B_1\cup \overline{B_1}$.
 Since $y'\ne y_s$, we have $|B_1|=|N_G(y')|-|\CC_{1y'}|\ge |N_G(y')\cap V_{\Delta}\cap (S\cup V(\CC_{\ge 3}))|$.
So
$$|B_1|\ge |N_G(y')\cap V_{\Delta}\cap (S\cup V(\CC_{\ge 3}))|\ge
\left\{
  \begin{array}{ll}
    |S_0|+(m_{11}+1)/2+m_{12}+m_2+m_3, &  \hbox{ if $y'\in T_1^1$\,(Claim~\ref{c111});} \\
    |S_0|+(m_{11}+1)/2+m_{12}+m_2+m_3-1, &  \hbox{ if $y'\in T_1^2$\,(Claim~\ref{c112});} \\
    |S_0|+m_{11}/2+m_{12}+m_2+m_3-1, & \hbox{ if $y'\in T_2$\,(Claim~\ref{c12}).}\\
    \end{array}
\right.
$$

In notching that  if $m_2>0$ then $m_2\ge 2$, using the above lower bounds on $|B_1|$,
we claim the following.
    \begin{equation}\label{B1}
    |B_1|\ge
    \left\{
  \begin{array}{ll}
    |S_0|+(m_{11}+1)/2+m_3, &  \hbox{(a) if $m_{12}=0,1$ and $ m_2=0$;} \\
    |S_0|+m_{11}/2+m_2+m_3-1, &  \hbox{(b) if $m_{12}=0, m_2\ge 2$;} \\
    |S_0|+(m_{11}+1)/2+m_{12}+m_3-1, & \hbox{(c) if $m_{12}\ge 2$ and $m_2=0$;}\\
    |S_0|+m_{11}/2+m_{12}+m_2+m_3-1, & \hbox{(d) if $m_{12}\ge 1$ and $m_2\ge 2$}.
  \end{array}
\right.
    \end{equation}
We verify (a). Notice that  when $m_{12}=0$ and $m_2=0$, $T_1^2=T_2=\emptyset$, which implies $y'\in T_1^1$.
By Claim~\ref{c111} we get $|B_1|\ge |S_0|+(m_{11}+1)/2+m_3$.
When $m_{12}=1$ and $m_2=0$, $y'\in T_1^1\cup T_1^2$. Then $|B_1|\ge \min\{|S_0|+(m_{11}+1)/2+1+m_3, |S_0|+(m_{11}+1)/2+1+m_3-1\}=|S_0|+(m_{11}+1)/2+m_3$.
Similarly, we can verify (b), (c) and (d).

By Claim~\ref{Hdegree} and~(\ref{degreeHrefine}), for each edge $xy\in E(H')$ with
$x\in B$ and $y\in A$, we have three cases:
\begin{itemize}
  \item [(i)]$d_H(y)+|\CC_{1y}|\ge 3\ge d_H(x)$ if $x\in B_1\cap U^{\CC_{\ge 3}^2}$,
where $|B_1\cap U^{\CC_{\ge 3}^2}|=|N_H(y')\cap U^{\CC_{\ge 3}^2}|=|N_G(y')\cap V(\CC_{\ge 3}^2)|\le |\CC_{\ge 3}^2|=m_3$ by (3) of
 Lemma~\ref{minimal_barrier};
  \item [(ii)] $d_H(y)+|\CC_{1y}|\ge d_H(x)+1$ if $x\in  B_1\cap S_0$; and
  \item [(iii)] $d_H(y)+|\CC_{1y}|\ge d_H(x)+2$
 if $x\in B_1-U^{\CC_{\ge 3}^2}-S_0$.
\end{itemize}
As $|A|>|B|$ and $M$ saturates $B$, $A-V(M)\ne \emptyset$.
Let $y^*\in A-V(M)$. Then  $d_H(y^*)+|\CC_{1y^*}|\ge 3$ by~(\ref{degreeHrefine}).
Hence
\begin{eqnarray}
&& e_H(A,B) = \sum\limits_{xy\in M\atop x\in B, y\in A} d_{H'}(x)\le \sum\limits_{xy\in M\atop x\in B, y\in A} d_{H}(x) \label{eab1} \\
   &=& \sum\limits_{xy\in M\atop x\in B_1, y\in A} d_H(x) +\sum\limits_{xy\in M\atop x\in \overline{B_1}, y\in A} d_H(x)\nonumber \\
   &\le & \sum\limits_{xy\in M\atop x\in B_1\cap S_0, y\in A} (d_H(y)+|\CC_{1y}|-1)+\sum\limits_{xy\in M\atop x\in B_1\cap U^{\CC^2_{\ge3}}, y\in A} (d_H(y)+|\CC_{1y}|)+ \nonumber\\
   & & \sum\limits_{xy\in M\atop x\in B_1-(S_0\cup U^{\CC^2_{\ge3}}), y\in A} (d_H(y)+|\CC_{1y}| -2)+\sum\limits_{xy\in M\atop x\in \overline{B_1}, y\in A} (d_H(y)+|\CC_{1y}|)\nonumber\\
 &\le &\sum\limits_{xy\in M\atop x\in B, y\in A} d_H(y)-|B_1\cap S_0|-2|B_1-(S_0\cup U^{\CC^2_{\ge3}})|+\sum\limits_{ y\in A\cap V(M)}|\CC_{1y}|\nonumber\\
 &< &\sum\limits_{xy\in M\atop x\in B, y\in A} d_H(y)-|B_1\cap S_0|-2|B_1-(S_0\cup U^{\CC^2_{\ge3}})|+\sum\limits_{ y\in A\cap V(M)\cup \{y^*\}}|\CC_{1y}|+d_H(y^*). \label{eab}
\end{eqnarray}
As $|B_1\cap S_0|\le |S_0|$, and $|B_1\cap U^{\CC^2_{\ge3}}|\le m_3$,
$2|B_1-(S_0\cup U^{\CC^2_{\ge3}})|\ge 2(|B_1|-|B_1\cap S_0|-|B_1\cap U^{\CC^2_{\ge3}}|)\ge 2(|B_1|-|S_0|-m_3)$. Hence, by~(\ref{B1})
$$2|B_1-(S_0\cup U^{\CC^2_{\ge3}})|\ge
\left\{
  \begin{array}{ll}
    m_{11}+1, & \hbox{(a) if $m_{12}=0,1, m_2=0$;} \\
    m_{11}+2m_2-2\ge m_{11}+m_2, & \hbox{(b) if $m_{12}=0, m_2\ge 2$;} \\
    m_{11}+2m_{12}-1\ge m_{11}+m_{12}, & \hbox{(c) if $m_{12}\ge 2$ and $m_2=0$;}\\
   m_{11}+2m_{12}+2m_2-2\ge m_{11}+m_{12}+m_2, & \hbox{(d) if $m_{12}\ge 1$ and $m_2\ge 2$}.
  \end{array}
\right.
$$
On the other hand,
$$\sum\limits_{ y\in A\cap V(M)\cup \{y^*\}}|\CC_{1y}|\le |\CC_1|=|\CC^1_{11}|+|\CC^2_{11}|+|\CC_{12}|\le
\left\{
  \begin{array}{ll}
    m_{11}+1, &  \hbox{(a) if $m_{12}=0,1$ and $ m_2=0$;} \\
    m_{11}+m_2, &  \hbox{(b) if $m_{12}=0$;} \\
    m_{11}+m_{12}, & \hbox{(c) if $m_2=0$;}\\
   m_{11}+m_{12}+m_{2}, & \hbox{(d) otherwise}.
  \end{array}
\right.
$$
So \,\, $-|B_1\cap S_0|-2|B_1-(S_0\cup U^{\CC^2_{\ge3}})|+\sum\limits_{ y\in A\cap V(M)\cup \{y^*\}}|\CC_{1y}|\le 0$,
and thus
from  inequalities~(\ref{eab1}) and (\ref{eab}), we get
\begin{eqnarray*}
e_H(A,B) & \le  & \sum\limits_{xy\in M\atop x\in B, y\in A} d_H(x) <  \sum\limits_{xy\in M\atop x\in B, y\in A} d_H(y)+d_{H}(y^*)\le  e_H(A,B),
\end{eqnarray*}
achieving  a contradiction.
\qed

The proof of Theorem~\ref{main} is then completed.
\qqed

\textbf{{\noindent \large Acknowledgements}}

The authors would like to  thank  Rong Luo  and Yue Zhao
for their discussions.

\bibliographystyle{plain}
\bibliography{SSL-BIB}

\end{document}